\documentclass[12pt]{article}

\newtheorem{defi}{Definition}

\newtheorem{prop}[defi]{Proposition}
\newtheorem{thm}[defi]{Theorem}

\usepackage{amsmath,amssymb}
\usepackage[title]{appendix}
\usepackage{geometry}
\title{Note on E-Polynomials Associated to $\mathbb{Z}_4$-codes}

\author{Nur Hamid}

\date{}

\begin{document}

\maketitle

\begin{abstract}
The invariant theory of finite groups can connect the coding theory to the number theory. 
In this paper, 
under this conformity, 
we obtain the minimal generators of the rings of E-polynomials constructed from the groups related to $\mathbb{Z}_4$-codes. 
In addition, 
we determine the generators of the invariant rings appearing by E-Polynomials and complete weight enumerators of Type II $\mathbb{Z}_4$-codes. 
\end{abstract}

\section{Introduction}

Our study is inspired by the idea of Motomura and Oura \cite{MotomuraOura}.
In their paper, 
they introduced the E-polynomials associated to the $\mathbb{Z}_4$-codes. 
They determined both the ring and the field structures generated by that E-polynomials. 
E-polynomials associated to the binary codes were investigated in a previously conducted study ({\it see} \cite{OuraM}).
In the present paper, 
we deal with $\mathbb{Z}_4$-codes. 
Then, we define an E-polynomial with respect to the complete weight enumerator of $\mathbb{Z}_4$-codes and 
show that the ring generated by them is minimally generated by E-polynomials of the following weights: 
\[
8, 16, 24, 32, 40, 48, 56, 64, 72, 80.
\]
It seems that the ring generated by E-polynomials is not sufficient to generate the invariant ring for the finite group $G^8$ defined in the next section. 
By combining the E-polynomials and the complete weight enumerators of $\mathbb{Z}_4$-codes,
we present the generators of that invariant ring.

We denote by $\mathbb{C}$ the field of complex number as usual.
Let $A_w$ be a finite-dimensional vector space over $\mathbb{C}$.
We write the dimension formula of $A$ by the formal series
\[
\sum_{w=0}^{\infty} (\dim A_w) t^w.
\]
For the dimension formulas and the basic theory of E-polynomials used herein,
we refer to references \cite{BannaiEtAl} and \cite{MotomuraOura}.
For the computations, we use Magma \cite{bosmaetal} and SageMath \cite{sagemath}.

\section{Preliminaries}
We denote a primitive 8-th root of unity by $\eta_8$. 
Following the notation used in \cite{BannaiEtAl}, let $G$ be a finite matrix group generated by 
\[
\frac{\eta_8}{2}
\left(
\begin{array}{rrrr}
1 &  1 &  1 &  1 \\
1 &  i & -1 & -i \\
1 & -1 &  1 & -1 \\
1 & -i & -1 &  i 
\end{array}
\right)
\]
and $\text{diag} \left[1,\eta_8,-1,\eta_8,\right]$. 
Let $G^8$ be a matrix group generated by $G$ and $\text{diag}\left[ \eta_8,\eta_8,\eta_8,\eta_8 \right]$. 
The group $G$ is of order 384, whereas $G^8$ is of order 1536.
We denote by $\mathfrak{R}$ and $\mathfrak{R}^8$ the invariant rings of $G$ and $G^8$, respectively:
\[
\mathfrak{R}=\mathbf{C}[t_0,t_1,t_2,t_3]^{G},
\]
\[
\mathfrak{R}^8=\mathbf{C}[t_0,t_1,t_2,t_3]^{G^8}
\]
under an action of such matrices on the polynomial ring of four variables $t_0$, $t_1$, $t_2$, and $t_3$.
The dimension formulas of $\mathfrak{R}$ and $\mathfrak{R}^8$ are given as follows:
\[
\sum_w \left(\text{dim} \mathfrak{R}_w \right) t^w = 
\frac{\left( 1+t^2\right) \left( 1+t^{4}\right) \left(1-t^2+2t^8+2t^{16}-t^{18}+t^{24}\right)}
{\left( 1-t^8\right)^3 \left( 1-t^{24} \right)},
\]
\[
\sum_w \left(\text{dim} \mathfrak{R}_w^8 \right) t^w = 
\frac{\left( 1+t^8\right) \left( 1+t^{16}\right)^2}{\left( 1-t^8\right)^3 \left( 1-t^{24} \right)}.
\]
In the next section, we present a fundamental theory of code that can help us obtain the generators of ring $\mathfrak{R}^8$.

\section{Codes} \label{secCode}

Let $ \mathbb{Z}_4$ be the ring of integers modulo 4.
A code $C$ over $ \mathbb{Z}_4 $ of length $n$, called $\mathbb{Z}_4$-code,
is an additive subgroup of $\mathbb{Z}_4^n$.
The inner product of two elements $a,b \in C$ on $\mathbb{Z}_4^n$ is given by
\[ (a,b)=a_1b_1 + a_2 b_2 + ... + a_n b_n  \mod 4 
 \] 
where $a=(a_1,a_2,...,a_n)$ and $b=(b_1,b_2,...,b_n)$.
The dual of $C$ is code $C^\perp$ satisfying
\[ 
C^\perp = \lbrace y \in \mathbb{Z}_4^n | (x,y) \equiv 0 \mod 4, \forall x \in C \rbrace.
\]
We say that $C$ is self-orthogonal if $C \subset C^\perp$ and self-dual if $C = C^\perp$.
A code $C$ is called \textit{Type II} if it is self-dual and satisfies
\[ 
(x,x) \equiv 0 \mod 8
 \]
for all $x \in C$. Type II $\mathbb{Z}_4$-code can only exist when its length is multiple of 8.
 
There are several types of weight enumerators associated with a $\mathbb{Z}_4$-code.
In this paper, we deal with complete weight enumerators. 

The \textit{complete weight enumerator} (CW) of a $\mathbb{Z}_4$-code $C$ is defined by
 \[ 
 CW_C (t_0, t_1, t_2, t_3)=\sum_{c \in C} t_0^{n_0(c)} t_1^{n_1(c)} t_2^{n_2(c)} t_3^{n_4(c)} 
  \]
where $n_i (c)$ denotes the number of $c$ components which are equivalent to $i$ modulo 4.
For every Type II $\mathbb{Z}_4$-code, $CW_C(t_0,t_1,t_2,t_3)$ is $G^8$-invariant. 
From the dimension formula, 
we have the following proposition.

\begin{prop} \label{pr_atmostgenerator}
The invariant ring $\mathfrak{R}^8$ can be generated by the set of complete weight enumerators of Type II $\mathbb{Z}_4$-codes consisting of at most 
\begin{center}
4 codes of length 8,\\
2 codes of length 16,\\
3 codes of length 24, \\
1 code of length 32, \\
1 code of length 40. 
\end{center}
\end{prop}

We denote by $p_{8a}, p_{8b}, o_8, k_8, p_{16a}, p_{16b}, q_{24a}, q_{24b}, g_{24}, q_{32}$ the complete weight enumerators of some codes. 
The numbers written as subscript denote the weight of each polynomial.
The codes $o_8$, $k_8$, and $g_{24}$ are  
known as octacode,  Klemm code, and Golay code, respectively.
The generator matrices of the complete weight enumerators which are denoted by $p$ are taken from \cite{PlessEtAl}.
We give the generator matrices of other complete weight enumerators in Appendix \ref{appendMatGen}.
The following are the explicit forms of some complete weight enumerators: 
\begin{align*}
p_{8a}  = \ &t_{0}^{8} + 4 t_{0}^{3} t_{1}^{4} t_{2} + 12 t_{0}^{6} t_{2}^{2} + 4 t_{0} t_{1}^{4} t_{2}^{3} + 38 t_{0}^{4} t_{2}^{4} + 12 t_{0}^{2} t_{2}^{6} + t_{2}^{8} + 4 t_{1}^{7} t_{3} + 16 t_{0}^{3} t_{1}^{3} t_{2} t_{3} \\
& + 16 t_{0} t_{1}^{3} t_{2}^{3} t_{3} + 24 t_{0}^{3} t_{1}^{2} t_{2} t_{3}^{2} + 24 t_{0} t_{1}^{2} t_{2}^{3} t_{3}^{2} + 28 t_{1}^{5} t_{3}^{3} + 16 t_{0}^{3} t_{1} t_{2} t_{3}^{3} + 16 t_{0} t_{1} t_{2}^{3} t_{3}^{3} \\
& + 4 t_{0}^{3} t_{2} t_{3}^{4} + 4 t_{0} t_{2}^{3} t_{3}^{4}  + 28 t_{1}^{3} t_{3}^{5} + 4 t_{1} t_{3}^{7}, \\
p_{8b} = \ & t_{0}^{8} + 8 t_{0}^{3} t_{1}^{4} t_{2} + 12 t_{0}^{6} t_{2}^{2} + 8 t_{0} t_{1}^{4} t_{2}^{3} + 38 t_{0}^{4} t_{2}^{4} + 12 t_{0}^{2} t_{2}^{6} + t_{2}^{8} + 16 t_{1}^{6} t_{3}^{2} + 48 t_{0}^{3} t_{1}^{2} t_{2} t_{3}^{2} \\
& + 48 t_{0} t_{1}^{2} t_{2}^{3} t_{3}^{2} + 32 t_{1}^{4} t_{3}^{4} + 8 t_{0}^{3} t_{2} t_{3}^{4} + 8 t_{0} t_{2}^{3} t_{3}^{4} + 16 t_{1}^{2} t_{3}^{6},\\
k_8  = \ & t_{0}^{8} + t_{1}^{8} + 28 t_{0}^{6} t_{2}^{2} + 70 t_{0}^{4} t_{2}^{4} + 28 t_{0}^{2} t_{2}^{6} + t_{2}^{8} + 28 t_{1}^{6} t_{3}^{2} + 70 t_{1}^{4} t_{3}^{4} + 28 t_{1}^{2} t_{3}^{6} + t_{3}^{8}, \\
o_8  = \ & t_{0}^{8} + t_{1}^{8} + 14 t_{0}^{4} t_{2}^{4} + t_{2}^{8} + 56 t_{0}^{3} t_{1}^{3} t_{2} t_{3} + 56 t_{0} t_{1}^{3} t_{2}^{3} t_{3} + 56 t_{0}^{3} t_{1} t_{2} t_{3}^{3} + 56 t_{0} t_{1} t_{2}^{3} t_{3}^{3} + 14 t_{1}^{4} t_{3}^{4} + t_{3}^{8},\\
p_{16a} = \ & t_{0}^{16} + 30 t_{0}^{8} t_{1}^{8} + t_{1}^{16} + 140 t_{0}^{12} t_{2}^{4} + 420 t_{0}^{4} t_{1}^{8} t_{2}^{4} + 448 t_{0}^{10} t_{2}^{6} + 870 t_{0}^{8} t_{2}^{8} + 30 t_{1}^{8} t_{2}^{8} \\
&+ 448 t_{0}^{6} t_{2}^{10} + 140 t_{0}^{4} t_{2}^{12} + t_{2}^{16} + 3360 t_{0}^{6} t_{1}^{6} t_{2}^{2} t_{3}^{2} + 6720 t_{0}^{4} t_{1}^{6} t_{2}^{4} t_{3}^{2} + 3360 t_{0}^{2} t_{1}^{6} t_{2}^{6} t_{3}^{2} \\
&+ 420 t_{0}^{8} t_{1}^{4} t_{3}^{4}  + 140 t_{1}^{12} t_{3}^{4} + 6720 t_{0}^{6} t_{1}^{4} t_{2}^{2} t_{3}^{4} + 19320 t_{0}^{4} t_{1}^{4} t_{2}^{4} t_{3}^{4} + 6720 t_{0}^{2} t_{1}^{4} t_{2}^{6} t_{3}^{4} \\
& + 420 t_{1}^{4} t_{2}^{8} t_{3}^{4}  + 448 t_{1}^{10} t_{3}^{6} + 3360 t_{0}^{6} t_{1}^{2} t_{2}^{2} t_{3}^{6} + 6720 t_{0}^{4} t_{1}^{2} t_{2}^{4} t_{3}^{6} + 3360 t_{0}^{2} t_{1}^{2} t_{2}^{6} t_{3}^{6} \\
&+ 30 t_{0}^{8} t_{3}^{8}  + 870 t_{1}^{8} t_{3}^{8} + 420 t_{0}^{4} t_{2}^{4} t_{3}^{8} + 30 t_{2}^{8} t_{3}^{8} + 448 t_{1}^{6} t_{3}^{10} + 140 t_{1}^{4} t_{3}^{12} + t_{3}^{16}.
\end{align*}
Since other weight enumerators are too large, 
we do not write them.

Let $\mathfrak{W}$ be a ring generated by the complete weight enumerators aforementioned:
\[
\mathfrak{W} = \mathbb{C}[p_{8a}, p_{8b}, o_8, k_8, p_{16a}, p_{16b}, q_{24a}, q_{24b}, g_{24}, q_{32}].
\]
By obtaining the dimension of $\mathfrak{W}$,
we have the following result. 
\begin{thm} \label{th_generatedenum}
The invariant ring $\mathfrak{R}^8$ can be generated by $\mathfrak{W}$.
\end{thm}

\textit{Proof.}
By Proposition \ref{pr_atmostgenerator},
we generate $\mathfrak{W}$ by utilizing some complete weight enumerators of non-equivalent codes. 
Then, 
we compute the dimension of $\mathfrak{W}$.
The dimension of each $\mathfrak{W}_k$ is shown in Table \ref{tb_dimeEnum}. 
This completes the proof of Theorem \ref{th_generatedenum}.

\begin{table} [ht]
\centering
\caption{The dimensions of $\mathfrak{R}^8_k$ and $\mathfrak{W}_k$}
\begin{tabular}{c|c|c|c|c|c}
$k$ & 8& 16 & 24 & 32 & 40 \\
\hline
dim$\mathfrak{R}^8_k$ & 4 & 11 & 25 & 48 & 83 \\
dim$\mathfrak{W}$ & 4 & 11 & 25 & 48 & 83  
\end{tabular}
\label{tb_dimeEnum}
\end{table}

It is noteworthy that we do not need to use the code of length 40.
On the next section, 
we shall give the generators of $\mathfrak{R}^8$ by the weight enumerators of Type II $\mathbb{Z}_4$-codes and E-polynomials. 

\section{E-Polynomials}

Let $\mathbf{t}$ be a column vector that comprises the following: $t_0$, $t_1$, $t_2$, and $t_3$.
An E-polynomial of weight $k$ for $G$ is defined by
\[
\varphi_k^G = \varphi_k^G(\mathbf{t})=\frac{1}{|G|} \sum_{\sigma \in G}(\sigma_0 \mathbf{t})^k
=\frac{|K|}{|G|} \sum_{ K\backslash G \ni \sigma}(\sigma_0 \mathbf{t})^k
\]
where
\[
K= \lbrace \begin{pmatrix}
1 & 0 & 0 & 0\\
\star & \star & \star & \star\\
\star & \star & \star & \star \\
\star & \star & \star & \star
\end{pmatrix} \in G
\rbrace
\]
and $\sigma_0$ is the first row of $\sigma$.
We apply the same definition for $G^8$.
The subgroup $K$ of $G$ is of order 8 and $K$ of $G^8$ is of order 16.
For simplicity, we denote by $\varphi_k$ without specifying the group. 
We denote by $\mathfrak{E}$ and $\mathfrak{E}^8$ the rings generated by $\varphi_k$s for the groups $G$ and $G^8$, respectively.

Denote by $\kappa$ the cardinality of $K \backslash G$. 
For clarity, 
we write $\kappa_G$ instead of $\kappa$ by including the group objected.
It is clear that $\kappa_G=48$ and $\kappa_{G^8}=96$.

\begin{thm} \label{finitelygenerated}
(1) 
The ring $\mathfrak{E}$ is generated by $\varphi_k$ where
\[
k \equiv 0 \mod 4, \quad 8 \leq k \leq 48.
\]
(2)
The ring $\mathfrak{E}^8$ is generated by $\varphi_k$ where
\[
k \equiv 0 \mod 8, \quad 8 \leq k \leq 96.
\]
\end{thm}

\textit{Proof.} 
(1) For each representative $\sigma_i $ of $ K \backslash G$ ($1 \leq i \leq \kappa$), let $x_i=\sigma_i' \mathbf{t}$, where $\sigma_i'$ is the first row of $\sigma_i$.
Then, every $\varphi_i$ can be expressed in $\mathbb{C}[x_1,\ldots,x_\kappa]$.
By the fundamental theorem of symmetric polynomials, 
every $\varphi_i$ can be written uniquely in $\varepsilon_i, \ldots, \varepsilon_\kappa \in \mathbb{C}[x_1, \ldots, x_\kappa]$
where
\[
\varepsilon_r = \sum_{i_1<i_2<\ldots<i_r} x_{i_1} x_{i_2} \ldots x_{i_r}, \quad (1\leq r \leq \kappa).
\]
We mention that $\varphi_4$=0.
This completes the proof. 

(2) The proof follows similarly that of Theorem \ref{finitelygenerated} (1).

Theorem \ref{finitelygenerated} informs us that the rings $\mathfrak{E}$ and $\mathfrak{E^8}$ are finitely generated. 
Hence, we can find their minimal generators. 
In the next theorem, 
we determine the generators of both $\mathfrak{E}$ and $\mathfrak{E^8}$.

\begin{thm} \label{generator}
(1) $\mathfrak{E}$ is minimally generated by the E-polynomials of weights
\[
8,12,16,20,24,28,32,40,48.
\]
(2) $\mathfrak{E}^8$ is minimally generated by the E-polynomials of weights
\[
8,16,24,32,40,48,56,64,72,80.
\]
\end{thm}

\textit{Proof.} 
For each $k$, we construct the rings $\mathfrak{E}_k$ and $\mathfrak{E}_k^8$.
Then, we determine whether $\varphi_k$ is generator or not. 
The dimensions of each $\mathfrak{E}$ and $\mathfrak{E}^8$ are demonstrated in Tables \ref{tb_dime} and \ref{tb_dime8}.
This completes the proof of Theorem \ref{generator}.

\begin{table} 
\centering
\caption{The dimensions of $\mathfrak{R}_k$ and $\mathfrak{E}_k$}
\begin{tabular}{c|c|c|c|c|c|c|c|c|c|c|c}
$k$ & 8&12 & 16 & 20 & 24 & 28 & 32 & 36 & 40 & 44 & 48 \\
\hline
dim$\mathfrak{R}_k$ & 4 & 3 & 16 & 11 & 25 & 27 & 48 & 54 & 83 & 94 & 133 \\
dim$\mathfrak{E}_k$ & 1 & 1 & 2 & 2 & 4 & 4 & 4 & 7 & 7 &  10 & 18
\end{tabular}
\label{tb_dime}
\end{table}

\begin{table} 
\centering
\caption{The dimensions of $\mathfrak{R}^8_k$ and $\mathfrak{E}^8_k$}
\begin{tabular}{c|c|c|c|c|c|c|c|c|c|c|c|c}
$k$ & 8& 16 & 24 & 32 & 40 & 48 & 56 & 64 & 72 & 80 & 88 & 96\\
\hline
dim$\mathfrak{R}^8_k$ & 4 & 11 & 25 & 48 & 83 & 133 & 200 & 287 & 397 & 532 & 695 & 889 \\
dim$\mathfrak{E}^8_k$ & 1 & 2 & 3 & 5 & 7 & 11 & 15 & 22 & 30 &  42 & 52 & 61
\end{tabular}
\label{tb_dime8}
\end{table}

Now, we obtain the relation between $\mathfrak{E}^8$ and $\mathfrak{R}^8$.
If we look at Table \ref{tb_dime8},
the ring $\mathfrak{E}^8$ is not sufficient to generate $\mathfrak{R}^8$.
By combining $\mathfrak{R}^8$ and $\mathfrak{W}$,
we have the following theorem. 

\begin{thm} \label{thm_generators}
The invariant ring $\mathfrak{R}^8$ can be generated by $\mathfrak{E}^8$
and the complete weight enumerators 
\[
p_8,o_8,k_8,p_{16},p_{24},q_{24},p_{32}.
\]

More specifically, the set 
\[
\lbrace \varphi_k, p_8,o_8,k_8,p_{16},p_{24},q_{24},p_{32} \quad | \quad k=8,16,24 \rbrace
\]
generates ring $\mathfrak{R}^8$.
\end{thm}

\textit{Proof.} 
Denote by $\widetilde{\mathfrak{R}}$ the polynomial generated by $\mathfrak{E}^8$ and the complete weight enumerators aforementioned. 
Then we construct $\widetilde{\mathfrak{R}}_k$ for $k=8,16,\ldots,96$.
It follows that 
each $\varphi_k$ for $k\neq 8, 16, 24$ is linearly dependent. 
We compute the dimension of each $\widetilde{\mathfrak{R}}_k$ 
and write the results in Table \ref{tb_dimR8}.
This completes the proof.

\begin{table}[h] 
\centering
\caption{The dimensions of $\mathfrak{R}^8_k$ and $\widetilde{\mathfrak{R}}$}
\begin{tabular}{c|c|c|c|c|c}
$k$ & 8& 16 & 24 & 32 & 40 \\
\hline
dim$\mathfrak{R}^8_k$ & 4 & 11 & 25 & 48 & 83 \\
dim$\widetilde{\mathfrak{R}}$ & 4 & 11 & 25 & 48 & 83 \\
\end{tabular}
\label{tb_dimR8}
\end{table}

\textit{Acknowledgment.} The author would like to thank Prof. Manabu Oura for his advice and suggestions.   

\begin{appendices}
\section{Generator Matrices}
\label{appendMatGen}

The generator matrix of $q_{24a}$ is given by
\[
\left(\begin{array}{*{24}r}
1 & 0 & 1 & 0 & 1 & 1 & 1 & 0 & 0 & 1 & 1 & 0 & 0 & 0 & 2 & 1 & 0 & 0 & 1 & 0 & 1 & 1 & 0 & 1 \\
0 & 1 & 0 & 0 & 1 & 1 & 0 & 2 & 0 & 1 & 1 & 0 & 0 & 0 & 2 & 3 & 0 & 0 & 1 & 1 & 0 & 0 & 0 & 0 \\
0 & 0 & 2 & 0 & 0 & 0 & 0 & 0 & 0 & 0 & 0 & 0 & 0 & 0 & 0 & 2 & 0 & 0 & 0 & 2 & 0 & 0 & 2 & 0 \\
0 & 0 & 0 & 1 & 1 & 1 & 0 & 1 & 0 & 0 & 0 & 0 & 0 & 0 & 0 & 2 & 0 & 0 & 0 & 2 & 0 & 0 & 2 & 0 \\
0 & 0 & 0 & 0 & 2 & 0 & 0 & 2 & 0 & 0 & 0 & 0 & 0 & 0 & 0 & 0 & 0 & 0 & 0 & 2 & 0 & 0 & 0 & 2 \\
0 & 0 & 0 & 0 & 0 & 2 & 0 & 2 & 0 & 0 & 0 & 0 & 0 & 0 & 0 & 0 & 0 & 0 & 0 & 2 & 0 & 0 & 0 & 2 \\
0 & 0 & 0 & 0 & 0 & 0 & 2 & 0 & 0 & 0 & 0 & 0 & 0 & 0 & 0 & 2 & 0 & 0 & 0 & 2 & 0 & 0 & 2 & 0 \\
0 & 0 & 0 & 0 & 0 & 0 & 0 & 0 & 1 & 1 & 1 & 0 & 0 & 0 & 1 & 2 & 0 & 0 & 0 & 2 & 0 & 0 & 0 & 2 \\
0 & 0 & 0 & 0 & 0 & 0 & 0 & 0 & 0 & 2 & 0 & 0 & 0 & 0 & 2 & 0 & 0 & 0 & 0 & 2 & 0 & 0 & 0 & 2 \\
0 & 0 & 0 & 0 & 0 & 0 & 0 & 0 & 0 & 0 & 2 & 0 & 0 & 0 & 2 & 0 & 0 & 0 & 0 & 2 & 0 & 0 & 0 & 2 \\
0 & 0 & 0 & 0 & 0 & 0 & 0 & 0 & 0 & 0 & 0 & 1 & 1 & 1 & 2 & 1 & 0 & 0 & 0 & 1 & 1 & 1 & 2 & 1 \\
0 & 0 & 0 & 0 & 0 & 0 & 0 & 0 & 0 & 0 & 0 & 0 & 2 & 0 & 0 & 2 & 0 & 0 & 0 & 2 & 0 & 0 & 0 & 2 \\
0 & 0 & 0 & 0 & 0 & 0 & 0 & 0 & 0 & 0 & 0 & 0 & 0 & 2 & 0 & 2 & 0 & 0 & 0 & 2 & 0 & 0 & 0 & 2 \\
0 & 0 & 0 & 0 & 0 & 0 & 0 & 0 & 0 & 0 & 0 & 0 & 0 & 0 & 0 & 0 & 1 & 1 & 1 & 3 & 1 & 1 & 3 & 1 \\
0 & 0 & 0 & 0 & 0 & 0 & 0 & 0 & 0 & 0 & 0 & 0 & 0 & 0 & 0 & 0 & 0 & 2 & 0 & 0 & 0 & 0 & 2 & 0 \\
0 & 0 & 0 & 0 & 0 & 0 & 0 & 0 & 0 & 0 & 0 & 0 & 0 & 0 & 0 & 0 & 0 & 0 & 2 & 2 & 0 & 0 & 2 & 2 \\
0 & 0 & 0 & 0 & 0 & 0 & 0 & 0 & 0 & 0 & 0 & 0 & 0 & 0 & 0 & 0 & 0 & 0 & 0 & 0 & 2 & 0 & 0 & 2 \\
0 & 0 & 0 & 0 & 0 & 0 & 0 & 0 & 0 & 0 & 0 & 0 & 0 & 0 & 0 & 0 & 0 & 0 & 0 & 0 & 0 & 2 & 0 & 2
\end{array}\right)
\]

The generator matrix of $q_{24b}$ is given by
\[
\left(\begin{array}{*{24}r}
1 & 0 & 0 & 0 & 0 & 0 & 1 & 0 & 0 & 1 & 0 & 0 & 0 & 0 & 0 & 2 & 0 & 1 & 0 & 1 & 1 & 2 & 1 & 3 \\
0 & 1 & 1 & 0 & 0 & 0 & 0 & 2 & 0 & 1 & 0 & 0 & 0 & 0 & 0 & 2 & 0 & 1 & 0 & 1 & 1 & 0 & 1 & 1 \\
0 & 0 & 2 & 0 & 0 & 0 & 2 & 0 & 0 & 0 & 0 & 0 & 0 & 0 & 0 & 0 & 0 & 0 & 0 & 0 & 0 & 2 & 0 & 2 \\
0 & 0 & 0 & 1 & 1 & 1 & 2 & 1 & 0 & 0 & 0 & 0 & 0 & 0 & 0 & 0 & 0 & 0 & 0 & 0 & 0 & 2 & 0 & 2 \\
0 & 0 & 0 & 0 & 2 & 0 & 0 & 2 & 0 & 0 & 0 & 0 & 0 & 0 & 0 & 0 & 0 & 0 & 0 & 0 & 0 & 0 & 0 & 0 \\
0 & 0 & 0 & 0 & 0 & 2 & 0 & 2 & 0 & 0 & 0 & 0 & 0 & 0 & 0 & 0 & 0 & 0 & 0 & 0 & 0 & 0 & 0 & 0 \\
0 & 0 & 0 & 0 & 0 & 0 & 0 & 0 & 1 & 1 & 1 & 0 & 0 & 0 & 1 & 2 & 0 & 0 & 0 & 1 & 1 & 3 & 2 & 3 \\
0 & 0 & 0 & 0 & 0 & 0 & 0 & 0 & 0 & 2 & 0 & 0 & 0 & 0 & 0 & 0 & 0 & 0 & 0 & 0 & 0 & 0 & 0 & 2 \\
0 & 0 & 0 & 0 & 0 & 0 & 0 & 0 & 0 & 0 & 2 & 0 & 0 & 0 & 0 & 0 & 0 & 0 & 0 & 0 & 0 & 2 & 0 & 0 \\
0 & 0 & 0 & 0 & 0 & 0 & 0 & 0 & 0 & 0 & 0 & 1 & 1 & 1 & 0 & 1 & 0 & 0 & 0 & 0 & 0 & 0 & 0 & 2 \\
0 & 0 & 0 & 0 & 0 & 0 & 0 & 0 & 0 & 0 & 0 & 0 & 2 & 0 & 0 & 2 & 0 & 0 & 0 & 0 & 0 & 0 & 0 & 0 \\
0 & 0 & 0 & 0 & 0 & 0 & 0 & 0 & 0 & 0 & 0 & 0 & 0 & 2 & 0 & 2 & 0 & 0 & 0 & 0 & 0 & 0 & 0 & 0 \\
0 & 0 & 0 & 0 & 0 & 0 & 0 & 0 & 0 & 0 & 0 & 0 & 0 & 0 & 2 & 0 & 0 & 0 & 0 & 0 & 0 & 2 & 0 & 0 \\
0 & 0 & 0 & 0 & 0 & 0 & 0 & 0 & 0 & 0 & 0 & 0 & 0 & 0 & 0 & 0 & 1 & 1 & 1 & 0 & 0 & 0 & 1 & 2 \\
0 & 0 & 0 & 0 & 0 & 0 & 0 & 0 & 0 & 0 & 0 & 0 & 0 & 0 & 0 & 0 & 0 & 2 & 0 & 0 & 0 & 0 & 2 & 0 \\
0 & 0 & 0 & 0 & 0 & 0 & 0 & 0 & 0 & 0 & 0 & 0 & 0 & 0 & 0 & 0 & 0 & 0 & 2 & 0 & 0 & 2 & 2 & 2 \\
0 & 0 & 0 & 0 & 0 & 0 & 0 & 0 & 0 & 0 & 0 & 0 & 0 & 0 & 0 & 0 & 0 & 0 & 0 & 2 & 0 & 0 & 0 & 2 \\
0 & 0 & 0 & 0 & 0 & 0 & 0 & 0 & 0 & 0 & 0 & 0 & 0 & 0 & 0 & 0 & 0 & 0 & 0 & 0 & 2 & 0 & 0 & 2
\end{array}\right)
\]

The generator matrix of $q_{32}$ is given by
\[
\left( \begin{array}{r}
10101010011000000010001201012123\\
01001000011000000010001201001020\\
00200002000000000000000000000022\\
00011103000000000000000000013101\\
00002002000000000000000000002002\\
00000202000000000000000000000000\\
00000022000000000000000000000022\\
00000000111000120000000000002002\\
00000000020000200000000000002002\\
00000000002000200000000000002002\\
00000000000111210000000000000000\\
00000000000020020000000000000000\\
00000000000002020000000000000000\\
00000000000000001110001200002002\\
00000000000000000200002000002002\\
00000000000000000020002000000000\\
00000000000000000001112100000000\\
00000000000000000000200200000000\\
00000000000000000000020200000000\\
00000000000000000000000011111133\\
00000000000000000000000002002022\\
00000000000000000000000000200020\\
00000000000000000000000000020002\\
00000000000000000000000000000202
\end{array} \right)
\]

\section{Other E-polynomials}
Let $G$ and $H$ be the matrix groups described as follows:
\[
G=\langle \frac{1}{i \sqrt{3}}
\begin{pmatrix}
1 & 2 \\
1 & -1
\end{pmatrix},
\begin{pmatrix}
1 & 0 \\
0 & e^{\frac{2\pi i}{3}}
\end{pmatrix}
\rangle,
\]
\[
H=\langle
\begin{pmatrix}
1 & 2 & 2 \\
1 & \zeta + \zeta^4 & \zeta^2+\zeta^3\\
1 & \zeta^2+\zeta^3 & \zeta+\zeta^4
\end{pmatrix},
\begin{pmatrix}
1 & 0 & 0 \\
0 & \zeta^2 & 0 \\
0 & 0 & \zeta^3
\end{pmatrix},
- \begin{pmatrix}
1 & 0 & 0 \\
0 & 0 & 1 \\
0 & 1 & 0
\end{pmatrix}
\rangle.
\]

The group $G$ is of order 24, whereas $H$ is of order 120.
The group $G$ is related to the self-dual ternary codes, whereas $H$ is related to the ring of symmetric Hilbert modular form. 
The discussion on these group can be found in \cite{EbelingWolfgang}.

By utilizing the same method discussed, 
we have that the ring generated by E-polynomials $\varphi_k^G$s (respectively $\varphi_k^H$s) 
is minimally generated by E-polynomials $\varphi_4$ and $\varphi_6$ 
(respectively $\varphi_2$, $\varphi_6$, and $\varphi_{10}$).
Thus, we have that
\[
\mathfrak{E}(G) =\langle \varphi_4, \varphi_6 \rangle
\]
and
\[
\mathfrak{E}(H) =\langle \varphi_2, \varphi_6, \varphi_{10} \rangle.
\]

The following tables present the dimensions of $\mathfrak{E}$ for each group. 
\begin{table}[ht]
\centering
\caption{The dimensions of $\mathfrak{R}(G)_k$ and $\mathfrak{E}(G)_k$}
\begin{tabular}{c|c|c}
$k$ & 4 & 6 \\
\hline
dim$\mathfrak{R}_k$ & 1 & 1 \\
dim$\mathfrak{E}_k$ & 1 & 1 
\end{tabular}
\label{tabEpolyG}
\end{table}

\begin{table}[ht]
\centering
\caption{The dimensions of $\mathfrak{R}(H)_k$ and $\mathfrak{E}(H)_k$}
\begin{tabular}{c|c|c|c|c|c}
$k$ & 2 & 4 & 6 & 8 & 10 \\
\hline
dim$\mathfrak{R}^8_k$ & 1 & 1 & 2 & 2 & 3 \\
dim$\mathfrak{E}^8_k$ & 1 & 1 & 2 & 2 & 3 
\end{tabular}
\label{tabEpolyH}
\end{table}
From Tables \ref{tabEpolyG} and \ref{tabEpolyH}, 
we can conclude that $\mathfrak{E}(G)$ (respectively $\mathfrak{E}(H)$)
satisfies
\[
\dim \mathfrak{E}(G)_k = \dim \mathfrak{R}(G)_k
\]
\[
(\dim \mathfrak{E}(H)_l = \dim \mathfrak{R}(H)_l)
\]
for $k\geq 4$ and $k \equiv 0 \mod 2$ (respectively $l \equiv 0 \mod 2$).
The dimension formulas of $\mathfrak{E}(G)$ and $\mathfrak{E}(H)$ can be written as follows. 
\[
G : \quad \frac{1}{(1-t^4)(1-t^6)},
\]
\[
H : \quad \frac{1}{(1-t^2)(1-t^6)(1-t^{10})}.
\]

\end{appendices}

\end{document}